\documentclass{article}

\usepackage{spconf}

\usepackage{amsmath}
\usepackage{amssymb}
\usepackage{mathtools}

\usepackage{hyperref}

\usepackage[cal=boondox]{mathalfa}

\newcommand{\Reals}{\mathbb{R}}
\newcommand{\dint}{\mathrm{d}}
\newcommand{\Df}[2]{D_f\bigl(#1 \,\Vert\, #2\bigr)}
\newcommand{\DKL}[2]{D_\text{KL}\bigl(#1 \,\Vert\, #2\bigr)}

\makeatletter
\newcommand*{\rom}[1]{\expandafter\@slowromancap\romannumeral #1@}
\makeatother

\newcommand{\refeqq}[1]{\mathrel{\overset{\makebox[0pt]{\mbox{\normalfont\scriptsize\sffamily #1}}}{=}}}
\newcommand{\refpropto}[1]{\mathrel{\overset{\makebox[0pt]{\mbox{\normalfont\scriptsize\sffamily #1}}}{\propto}}}

\DeclareMathOperator{\tr}{tr}
\DeclareMathOperator{\mse}{mse}
\DeclareMathOperator{\mmse}{mmse}

\AtBeginDocument{
  \label{CorrectFirstPageLabel}
  
}

\long\def\symbolfootnote[#1]#2{\begingroup%
	\def\thefootnote{\fnsymbol{footnote}}\footnote[#1]{#2}\endgroup} 

\title{TIGHT MMSE BOUNDS FOR THE AGN CHANNEL \\ UNDER KL DIVERGENCE CONSTRAINTS ON THE INPUT DISTRIBUTION}
\twoauthors
 {Michael~Fau\ss, Abdelhak M.~Zoubir}
	{Technische Universit\"at Darmstadt\\
	Signal Processing Group\\
	64283 Darmstadt, Germany}
 {Alex~Dytso, H.~Vincent~Poor\sthanks{This work was supported in part by the U.S.~National Science Foundation under Grants CCF-1420575 and ECCS-1549881.}}
	{Princeton University\\
	Dept.~of Electrical Engineering\\
	Princeton, NJ 08544, USA}

\begin{document}

\ninept

\maketitle

\begin{abstract}
Tight bounds on the minimum mean square error for the additive Gaussian noise channel are derived, when the input distribution is constrained to be $\varepsilon$-close to a Gaussian reference distribution in terms of the Kullback--Leibler divergence. The distributions that attain the bounds are shown be Gaussian whose means are identical to that of the reference distribution and whose covariance matrices are defined implicitly via systems of matrix equations. The estimator that attains the upper bound is identified as a minimax optimal estimator that is robust against deviations from the assumed prior. The lower bound is shown to provide a potentially tighter alternative to the Cram\'{e}r--Rao bound. Both properties are illustrated with numerical examples.
\end{abstract}



\section{Introduction and Problem Formulation}
\label{sec:problem}

The mean square error (MSE) is a natural and commonly used measure for the accuracy of an estimator. The minimum  MSE (MMSE)  plays a central role in statistics \cite{Lehmann_1998, Dodge_2008}, information theory \cite{Guo_2011, Dytso_2018}, and signal processing \cite{Kay_1993, Dalton_2012, Dalton_2012b} and has been shown to have close connections to entropy and mutual information \cite{Guo_2005, Guo_2006}.

In this paper, lower and upper bounds on the MMSE are derived when the random variable of interest is contaminated by additive Gaussian noise and its distribution is constrained to be $\varepsilon$-close to a Gaussian reference distribution in terms of the Kullback--Leibler (KL) divergence. This problem is of interest in both information theory as well as robust statistics. More precisely, it is shown that the estimator that attains the upper bound is minimax robust in the sense that it minimizes the maximum MMSE over the set of feasible distributions. That is, within the specified  KL divergence ball, it is robust against arbitrary deviations of the prior from the nominal Gaussian case.  In addition, the lower bound provides a fundamental limit on the estimation accuracy and provides an alternative to the well-known Bayesian Cram\'{e}r--Rao bound. However, since it uses additional information about the KL divergence, it can be significantly tighter in some cases. 

More formally, let $(\mathbb{R}^K, \mathcal{B}^K)$ denote the $K$-dimensional Borel space. Consider an additive-noise channel
\begin{equation*}
  Y = X + N,
\end{equation*}
where $X$ and $N$ are independent $(\mathbb{R}^K, \mathcal{B}^K)$-valued random variables. Without loss of generality, $N$ is assumed to be zero-mean. For the purpose of this paper, it is useful to define the MSE as a function of an estimator $f$ and an input distribution $P_X$, i.e.,
\begin{equation*}
  \mse_{X | Y}(f, P_X) \coloneqq E_{P_{Y|X} P_X}\Bigl[ \left\lVert f(Y) - X \right\rVert^2\Bigr].
\end{equation*}
The MMSE is accordingly defined as
\begin{equation*}
  \mmse_{X | Y}(P_X) \coloneqq \inf_{f \in \mathcal{F}} \; \mse_{X | Y}(f, P_X),
\end{equation*}
where $\mathcal{F}$ denotes the set of all all feasible estimators, i.e.,
\begin{equation*}
  \mathcal{F} = \left\{ f\colon (\mathbb{R}^K, \mathcal{B}^K) \to (\mathbb{R}^K, \mathcal{B}^K) \right\}.
\end{equation*}
The problems investigated in this paper are
\begin{align}
  \label{eq:problem_sup}
  \sup_{P_X} \; \mmse_{X|Y}(P_X) \quad &\text{s.t.} \quad P_X \in \mathcal{P}_\varepsilon, \\
  \label{eq:problem_inf}
  \inf_{P_X} \; \mmse_{X|Y}(P_X) \quad &\text{s.t.} \quad P_X \in \mathcal{P}_\varepsilon.
\end{align}
 where the set of all feasible distribution is defined as 
\begin{equation*}
  \mathcal{P}_\varepsilon \coloneqq \bigl\{ P_X \mid \DKL{P_X}{P_0} \leq \varepsilon \bigr\}.
\end{equation*}
Note that $\mathcal{P}_\varepsilon$ is a KL ball centered at $P_0$ of radius $\epsilon$.  It is further assumed that
\begin{equation*}
  P_0 = \mathcal{N}(\mu_0,\Sigma_0),
\end{equation*}
where $\mathcal{N}(\mu,\Sigma)$ denotes the Gaussian distribution with mean $\mu$ and covariance $\Sigma$. All covariance matrices are assumed to be positive definite. 

The paper is organized as follows: the main result is stated in Section~\ref{sec:result}, followed by a brief discussion in Section~\ref{sec:discussion}. A proof is detailed in Section~\ref{sec:proof} and two illustrative numerical examples are presented in Section~\ref{sec:examples}.

A remark on notation: In what follows, $\Sigma_X$ and $\Sigma_N$ denote the covariance matrices of $X$ and $N$, respectively. Moreover, in order to keep the notation compact, the following matrices are introduced 
\begin{align*}
  W_X &\coloneqq \Sigma_X (\Sigma_X + \Sigma_N)^{-1}, \\
  D_X &\coloneqq \Sigma_N W_X^\text{T} W_X, \\
  W_N &\coloneqq \Sigma_N (\Sigma_X + \Sigma_N)^{-1}, \\
  D_N &\coloneqq \Sigma_X W_N^\text{T} W_N,
\end{align*}
where $W^\text{T}$ denotes the transpose of $W$. Note that $W_X + W_N = I$.

\section{Result}
\label{sec:result}

If $(\alpha^* ,\Sigma_X^*)$, with $\Sigma_X^*$ positive definite and $\alpha^* \geq 0$, solve
\begin{gather}
  \label{eq:opt_covariance}
  \Sigma_X = (I + \alpha D_N) \Sigma_0, \\
  \label{eq:opt_KL}
  \alpha \tr\bigl( D_N \bigr) - \log\det\bigl( I + \alpha D_N \bigr) = 2\varepsilon,
\end{gather}
then 
\begin{equation}
  \label{eq:px_opt_sup}
  P_X^* = \mathcal{N}\bigl(\mu_0, \Sigma_X^* \bigr)
\end{equation}
solves \eqref{eq:problem_sup}. Analogously, if $(\alpha^\dagger ,\Sigma_X^\dagger)$, with $\Sigma^\dagger$ positive definite and $\alpha^\dagger \leq 0$, solve \eqref{eq:opt_covariance} and \eqref{eq:opt_KL}, then 
\begin{equation}
  \label{eq:px_opt_inf}
  P_X^\dagger = \mathcal{N}\bigl(\mu_0, \Sigma_X^\dagger \bigr)
\end{equation}
solves \eqref{eq:problem_inf}. 

Since the optimal distributions are Gaussians of the form $P_X = \mathcal{N}\bigl(\mu_0, \Sigma_X \bigr)$, the MMSE estimators are in both cases given by
\begin{equation}
  \label{eq:MMSE_estimator}
  f(y) = W_X y + W_N \mu_0,
\end{equation}
where $y \in \Reals^K$ denotes the observed realization of $Y$. The corresponding MMSE calculates to
\begin{equation}
  \label{eq:MMSE}
  \mmse_{X|Y}(P_X) = \tr\bigl( D_X \bigr) + \tr\bigl( D_N \bigr).
\end{equation}
The lower and upper bound can be obtained by evaluating \eqref{eq:MMSE} at $\Sigma_X^*$ and $\Sigma_X^\dagger$, respectively.

\section{Discussion}
\label{sec:discussion}

Before proceeding with the proof of the presented bounds, some clarifying remarks are in order.

1) The result does not make a statement about the existence or the uniqueness of $\Sigma_X^*$ and $\Sigma_X^\dagger$. However, considering that the MMSE is a concave functional of the input distribution \cite[Corollary 1]{Wu_2012} and that the KL divergence is strictly convex in both arguments \cite[Theorem 2.7.2]{Cover_Thomas_2006}, we conjecture that the solution is guaranteed to exist. But, considering that the MMSE is not strictly concave \cite[Corollary 2]{Wu_2012}, it is likely not to be unique.

2) Solving \eqref{eq:opt_covariance} and \eqref{eq:opt_KL} for $\Sigma_X$ and $\alpha$ is non-trivial and, in general, requires the use of numerical solvers. For this purpose, it can be useful to rewrite \eqref{eq:opt_covariance} as
\begin{equation}
  \label{eq:opt_covariance_alternative}
  (I + \Sigma_N^{-1} \Sigma_X)^{\text{T}}(I + \Sigma_N^{-1} \Sigma_X)(I-\Sigma_0^{-1} \Sigma_X) + \alpha \Sigma_X = 0
\end{equation}
in order to avoid numerical instabilities that might arise from inverting $\Sigma_X$. A detailed derivation of \eqref{eq:opt_covariance_alternative} is omitted due to space constraints. Our limited experimental results suggest that both \eqref{eq:opt_covariance} and \eqref{eq:opt_covariance_alternative} can usually be solved using off-the-shelf algorithms. However, a detailed discussion is beyond the scope of this work.

3) A special case for which an analytical solution of \eqref{eq:opt_covariance} and \eqref{eq:opt_KL} can be given is a white signal in white noise, i.e., $\Sigma_0 = \sigma_0^2 I$ and $\Sigma_N = \sigma_N^2 I$. In this case it holds that
\begin{equation*}
  P_X^* = \mathcal{N}\bigl( \mu_0, s^* I \bigr) \quad \text{and} \quad P_X^\dagger = \mathcal{N}\bigl( \mu_0,  s^\dagger I \bigr) ,
\end{equation*}
where the variances $s^*$ and $s^\dagger$ are unique and given by
\begin{align}
  \label{eq:sigma_min}
  s^* &= -W_{0}\left(-e^{-\left(1+\frac{2}{K}\varepsilon\right)}\right) \sigma_0^2, \\
  \label{eq:sigma_max}
  s^\dagger &= -W_{-1}\left(-e^{-\left(1+\frac{2}{K}\varepsilon\right)}\right) \sigma_0^2.
\end{align}
Here $W_k$ denotes the $k$th branch of the Lambert W function \cite{Corless_1996}. The corresponding MMSE bounds calculate to
\begin{equation*}
  \frac{s^\dagger \sigma_N^2}{s^\dagger + \sigma_N^2} \leq \frac{\mmse_{X|Y}(P_X)}{K} \leq \frac{s^* \sigma_N^2}{s^* + \sigma_N^2} \quad \forall P_X \in \mathcal{P}_\varepsilon.
\end{equation*}
The proof of these results is straightforward and omitted for brevity.

4) The signal model can be extended to conditionally independent and identically distributed random variables $Y_1|X, \ldots, Y_n|X$ by letting
\begin{equation*}
  \Sigma_N \leftarrow \frac{1}{n} \Sigma_N \quad \text{and} \quad y \leftarrow \frac{1}{n}\sum_{i=1}^n y_i,
\end{equation*}
where $y_i$ denotes an observation of $Y_i|X$, $i=1,\ldots,n$.

\section{Proof}
\label{sec:proof}

The idea underlying the proof is to reformulate \eqref{eq:problem_sup} and \eqref{eq:problem_inf} as nested optimizations over the estimator $f$ and the distribution $P_X$. Only the solution of \eqref{eq:problem_sup} is presented in detail since the solution of \eqref{eq:problem_inf} can be given analogously.

Before detailing the proof of the main result, a result on the optimization of expected values under $f$-divergence constraints is derived. It simplifies the subsequent steps, but also constitutes a useful result in its own right.

\subsection{Bounding expectations under f-divergence constraints}

Consider the auxiliary problems
\begin{align}
  \label{eq:aux_problem_sup}
  \sup_{P_X} \; E_{P_X}\bigl[ h(X) \bigr] \quad \text{s.t.} \quad \Df{P_X}{P_0} \leq \varepsilon, \\
  \label{eq:aux_problem_inf}
  \inf_{P_X} \; E_{P_X}\bigl[ h(X) \bigr] \quad \text{s.t.} \quad \Df{P_X}{P_0} \leq \varepsilon,
\end{align}
where $h \colon (\mathbb{R}^K, \mathcal{B}^K) \to (\mathbb{R}, \mathcal{B})$ is a measurable function. Assuming $P_X$ and $P_0$ to be absolutely continuous with respect to a reference measure $\eta$, \eqref{eq:aux_problem_sup} and \eqref{eq:aux_problem_inf} can be reformulated in terms of the densities of $P_X$ and $P_0$ w.r.t.~$\eta$, namely,
\begin{align}
  \{ \sup_{p_X} \,,\, \inf_{p_X} \} \; \int h(x) p_X(x) \, \eta(\dint x) \\
  \text{s.t.} \quad \int_{\Reals^K} f\biggl(\frac{p_X(x)}{p_0(x)}\biggr) p_0(x) \, \eta(\dint x) &\leq \varepsilon \\ 
  \int_{\Reals^K} p_X(x) \, \eta(\dint x) &= 1 \\ 
  p_X(x) &\geq 0 \label{eq:nneg_constraint}  .
\end{align}
The equivalent unconstrained optimization problems are given by
\begin{equation}
  \label{eq:aux_problem_unconstrained}
  \{ \sup_{p_X} \,,\, \inf_{p_X} \} \; L(p_X; \lambda, \nu), 
\end{equation}
with
\begin{multline}
  L(p_X; \lambda, \nu)   \coloneqq   \\ \int_\Reals h(x) + \lambda \, f\biggl( \frac{p_X(x)}{p_0(x)} \biggr) p_0(x) + \nu \, p_X(x) \, \eta(\dint x), \label{eq:lagrange}
\end{multline}
and $\nu, \lambda \in \Reals$. Note that $\lambda \leq 0$ for \eqref{eq:aux_problem_sup} and $\lambda \geq 0$ for \eqref{eq:aux_problem_inf} and that the constraint \eqref{eq:nneg_constraint} is neglected since it turns out to be redundant. The Fr\'{e}chet derivative \cite{Bell_2014} of \eqref{eq:lagrange} w.r.t.~$p_X$ is given by \begin{equation*}
  L_{p_X}'(p_X; \lambda, \nu) = h(x) + \lambda f'\biggl(\frac{p_X(x)}{q(x)} \biggr) + \nu,
\end{equation*}
where $f'$ denotes the subderivative \cite[p.~36]{Clarke_1990} of $f$. A sufficient condition for $p_X^*$ to solve \eqref{eq:aux_problem_unconstrained} is that $L_{p_X}'(p_X; \lambda, \nu) = 0$, which yields
\begin{equation}
  \label{eq:aux_problem_opt}
  p_X^*(x) = q(x) g\biggl( -\frac{h(x)+\nu}{\lambda} \biggr) = q(x) g\big( \alpha h(x) + \beta \big), 
\end{equation}
where $\alpha = -1/\lambda$, $\beta = -\nu/\lambda$, and $g \colon \Reals \to \Reals$ is the generalized inverse \cite{Embrechts_2013} of $f'$, i.e.,
\begin{equation*}
  g(c) \coloneqq \inf \{ x \in \Reals \mid f'(x) \geq c \}.
\end{equation*}
Note that $\alpha \geq 0$ for \eqref{eq:aux_problem_sup} and $\alpha \leq 0$ for \eqref{eq:aux_problem_inf}. Finally, in order for $p_X^*(x)$ to solve the constrained problem \eqref{eq:aux_problem_sup} or \eqref{eq:aux_problem_inf}, $\alpha$ and $\beta$ need to be chosen such that $p_X^*(x)$ is a valid density and the $f$-divergence constraint is fulfilled with equality. The latter follows directly from the fact that the complementary slackness constraint of the Karush--Kuhn--Tucker conditions \cite{Brezhneva_2011} requires the $f$-divergence constraint to be satisfied with equality for all $\lambda \neq 0$ and hence for all $\alpha \in \Reals$. 

\subsection{Proof of the main result}

Consider the maximization in \eqref{eq:problem_sup} which can be written as the minimax problem
\begin{equation}
  \label{eq:nested_problem}
  \sup_{P_X \in \mathcal{P}_\varepsilon} \; \inf_{f \in \mathcal{F}} \; \mse_{X|Y}(f,P_X).
\end{equation}
A sufficient condition for $P_X^*$ and $f^*$ to solve \eqref{eq:nested_problem}, and hence \eqref{eq:problem_sup}, is that they satisfy the saddle point conditions \cite[Exercise 3.14]{Boyd_2004}
\begin{align}
  \label{eq:saddle_point_f} 
  \mse_{X | Y}(f^*,P_X^*) &\leq \mse_{X | Y}(f,P_X^*) & \forall f &\in \mathcal{F}, \\
  \label{eq:saddle_point_P}
  \mse_{X | Y}(f^*,P_X^*) &\geq \mse_{X | Y}(f^*,P_X) & \forall P_X &\in \mathcal{P}_\varepsilon.
\end{align}
The fact that $f$ in \eqref{eq:MMSE_estimator} minimizes the right hand side of \eqref{eq:saddle_point_f} follows directly from the definition of the MMSE \cite[Chapter 10.4]{Penny_2000}. In the remainder of the proof, it is shown that $P_X^*$ in \eqref{eq:px_opt_sup} satisfies \eqref{eq:saddle_point_P}.  

First, the right hand side of \eqref{eq:saddle_point_P} is written as
\begin{equation*}
  \mse_{X | Y}(f^*,P_X) = E_{P_X} \bigl[ h(X) \bigr],
\end{equation*}
where $h\colon \Reals^K \to \Reals$ is independent of $P_X$ and given by
\begin{align*}
  h(x) &= E_{P_{Y \mid X=x}} \Bigl[ \left\lVert W_X Y + W_N \mu_0 - x \right\rVert_2^2 \Bigr] \\
  &= E_{\mathcal{N}(x, \Sigma_N)} \Bigl[ \left\lVert W_X Y + W_N \mu_0 - x \right\rVert_2^2 \Bigr] \\
  &= \tr \bigl( \Sigma_N W_X^\text{T} W_X \bigr) + \left\lVert W_X x + W_N \mu_0 - x \right\rVert_2^2 \\
  &= \tr \bigl( D_X \bigr) + (x - \mu_0)^\text{T} W_N^\text{T} W_N (x-\mu_0).
\end{align*}
In order for $\DKL{P_X}{P_0}$ to be finite, $P_X$ needs to be absolutely continuous w.r.t.~$P_0$. Therefore, the problem
\begin{equation}
  \label{eq:opt_h}
  \sup_{P_X \in \mathcal{P}_\varepsilon}\; \mse_{X|Y}(f^*, P_X) = \sup_{P_X \in \mathcal{P}_\varepsilon} \; E_{P_X} \bigl[ h(X) \bigr]
\end{equation}
is of the form \eqref{eq:aux_problem_sup}, with $f(x) = x \log x$. Since the latter is strictly convex, the inverse function of its derivative is unique and is given by $g(c) = e^{c-1}$. Inserting $g(c)$ into \eqref{eq:aux_problem_opt} yields
\begin{align*}
  p_X^*(x) &= q(x) e^{\alpha h(x) + \beta -1} \\
  &\propto q(x) e^{\alpha (x - \mu_0)^\text{T} W_N^\text{T} W_N (x-\mu_0)} \\
  &\refpropto{(\,\rom{1}\,)} e^{-\frac{1}{2} (x-\mu_0)^\text{T} (\Sigma_0^{-1} - \alpha W_N^\text{T} W_N) (x-\mu_0)} \\
  &\propto e^{-\frac{1}{2}(x-\mu_0)^\text{T} \Sigma_X^{-1} (x-\mu_0)},
\end{align*}
where $\alpha \geq 0$ and
\begin{align}
  \label{eq:sigma_px}
  \Sigma_X^{-1} &= \Sigma_0^{-1} - \alpha W_N^\text{T} W_N. 
\end{align}
Note that, without loss of generality, $\alpha$ has been scaled by $1/2$ in (I). Multiplying \eqref{eq:sigma_px} by $\Sigma_X$ from the left, by $\Sigma_0$ from the right, and rearranging the terms yields the optimality condition in \eqref{eq:opt_covariance}. 

Knowing that $P_X$ and $P_0$ are Gaussians with identical means, the KL divergence $\DKL{P_X}{P_0}$ is given by 
\begin{align}
  \DKL{P_X}{P_0} &= \frac{1}{2}\left( \tr\bigl( \Sigma_X \Sigma_0^{-1} \bigr) - K - \log\det \Sigma_X \Sigma_0^{-1} \right) \notag \\
  &\refeqq{\eqref{eq:opt_covariance}} \frac{1}{2}\left( \alpha \tr\bigl( D_N \bigr) + \log\det\bigl( I + \alpha D_N \bigr) \right). \label{eq:gauss_kl}
\end{align}
Equating \eqref{eq:gauss_kl} with $\varepsilon$ yields the optimality condition \eqref{eq:opt_KL}. This concludes the proof.

The proof for the optimality of $P_X^\dagger$ follows analogously, the only difference being that in \eqref{eq:opt_h} the supremum is replaced by the infimum and, consequently, the sign of $\alpha$ is reversed.

\section{Numerical Examples}
\label{sec:examples}

In order to illustrate the usefulness of the bounds in Section~\ref{sec:result}, two examples are presented; the first in the context of robust MMSE estimation, the second in the context of estimation accuracy bounds. 

\subsection{Robust MMSE estimation}

In this example, the nominal MMSE estimator $f_0$ with $\Sigma_X = \Sigma_0$ and the robust MMSE estimator $f^*$  with $\Sigma_X = \Sigma_X^*$ are compared in terms of their performance under $P_0$ and under their respective least favorable distribution in the feasible set $\mathcal{P}_\varepsilon$. By definition, the least favorable distribution for the robust MMSE estimator $f^*$ is $P_X^*$. From \eqref{eq:sigma_px} it follows that the least favorable distribution for $f_0$ is given by a Gaussian with mean $\mu_0$ and covariance
\begin{equation*}
  \Sigma_X^{-1} = \Sigma_0^{-1} - \alpha (\Sigma_0 + \Sigma_N)^{-1} \Sigma_N^2 (\Sigma_0 + \Sigma_N)^{-1}, 
\end{equation*}
where $\alpha \geq 0$ needs to be chosen such that the KL divergence constraint is fulfilled with equality. In a slight abuse of notation, the MMSE of $f_0$ under the corresponding least favorable distribution is denoted by $\mmse_{X|Y}(f_0, P_X^*)$.

For the example, the matrix $\Sigma_0$ is assumed to be of size $K = 10$ and to admit a Toeplitz structure with entries
\begin{equation*}
  [\Sigma_0]_{ij} = e^{-0.9 \lvert i - j \rvert }, \quad i,j = 1, \ldots, 10.
\end{equation*} 
This implies $\tr(\Sigma_0) = K$. The noise is assumed to be white with variance $\sigma_N^2 = 1/\gamma$ so that $\gamma = \tr(\Sigma_0)/\tr(\Sigma_N)$ denotes the signal-to-noise ratio (SNR) under $P_0$. 

\begin{figure}[t]
  \centering
  \includegraphics{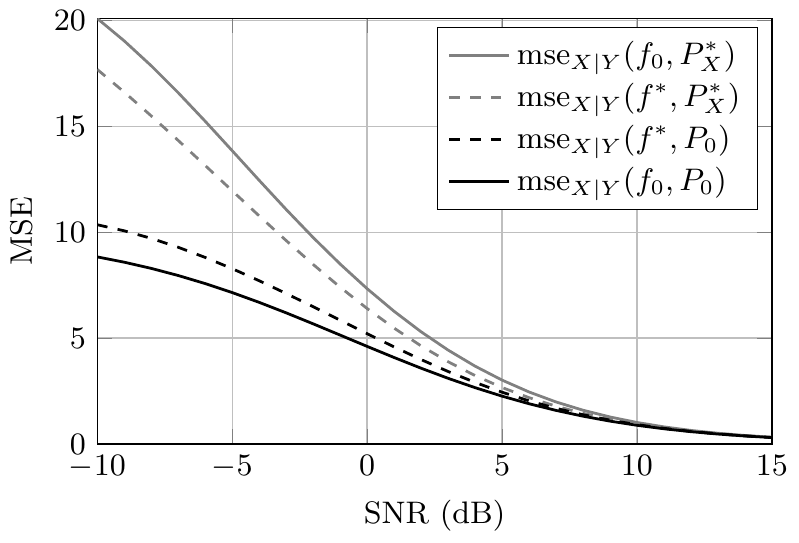}
  \caption{MSE vs.~SNR for the nominal MMSE estimator $f_0$ and the robust MMSE estimator $f^*$ under $P_0$ and the respective least favorable distribution $P_X^*$ for $\varepsilon = 2$.}
  \label{fig:MSE_vs_snr}
\end{figure}

In Fig.~\ref{fig:MSE_vs_snr}, the MSE of the estimators $f_0$ and $f^*$ under the nominal distribution $P_0$ and their respective least favorable distributions is plotted versus the SNR. The KL divergence tolerance is set to $\varepsilon = 2$. Especially in the low SNR regime, the robust estimator $f^*$ yields a gain of 1-3 \,dB in terms of the worst case MSE. The price for this improvement is a loss of comparable magnitude in terms of nominal performance.

\begin{figure}[t]
  \centering
  \includegraphics{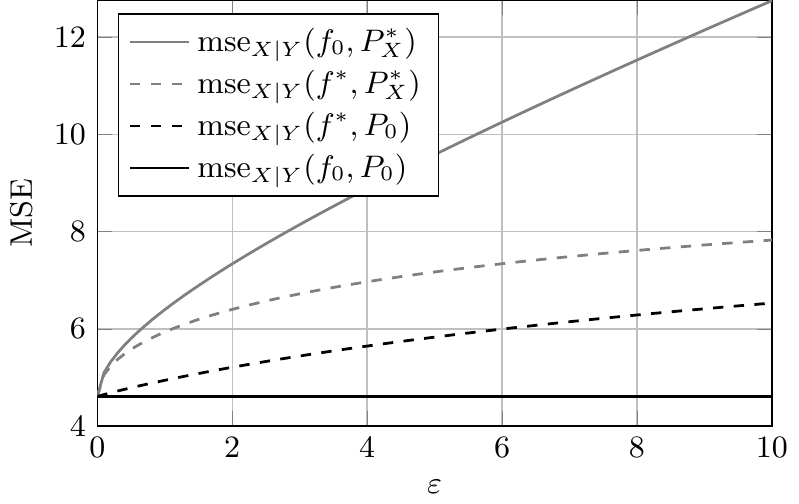}
  \caption{MSE vs.~$\varepsilon$ for the nominal MMSE estimator $f_0$ and the robust MMSE estimator $f^*$ under $P_0$ and the respective least favorable distribution $P_X^*$ at SNR $=$ 0\,dB.}
  \label{fig:MSE_vs_eps}
\end{figure}

The increased robustness of $f^*$ in comparison to $f_0$ becomes more apparent in Fig.~\ref{fig:MSE_vs_eps}, where the MSE is plotted versus the KL divergence tolerance $\varepsilon$ at an SNR of 0\,dB. It can clearly be seen how the worst case MSE scales differently for $f^*$ and $f_0$, which highlights the advantage of the robust estimator under large prior uncertainties. At the same time, the MSE of $f^*$ under the nominal distribution $P_0$ deteriorates with increasing $\varepsilon$. However, in this example, the gain in worst case MSE is significantly larger than the loss in nominal performance.

\subsection{Comparison to the Cram\'{e}r--Rao bound}

In this example, the usefulness of the lower bound in Section~\ref{sec:result} is illustrated by showing that it can provide a tighter alternative to the Cram\'{e}r--Rao lower bound. Let $K = 1$ and consider the zero-mean generalized Gaussian distribution $\mathcal{G}(a,p)$ with density function
\begin{equation*}
  \mathcal{g}(x \,|\,a, p) = \frac{p}{2a    \Gamma(1/p) } e^{-\left(\frac{\lvert x \rvert}{a}\right)^p},
\end{equation*}
where $\Gamma$ denotes the gamma function \cite{Davis_1959}, $a > 0$ is a scale parameter, and $p > 0$ determines the type of decay of the tails \cite{DytsoISIT2017}. Note that the generalized Gaussian reduces to a regular Gaussian for $p = 2$ and to a Laplace distribution for $p = 1$.   It is not difficult to verify that choosing $a= \sqrt{\frac{\Gamma(1/p)}{\Gamma(3/p)}} b$ implies that $E[ X^2]=b^2.$
  
Calculating the exact MMSE for $X \sim \mathcal{G}(a,p)$ is non-trivial in general so that lower bounds are of interest that are either analytical or can be calculated with a low computational effort. A commonly used bound on the MMSE is the Baysian Cram\'{e}r--Rao bound \cite{Gill_1995}, which  for the transformation $Y=\sqrt{\gamma}X+N$ with $N  \sim \mathcal{N}\bigl(0, 1 \bigr)$ is given by
\begin{equation*}
  \mmse_{X|Y}\bigl(\mathcal{G}(a,p)\bigr) \geq \frac{1}{\gamma + I\bigl( \mathcal{G}(a,p) \bigr)},
\end{equation*}
where
\begin{equation}
  \label{eq:gg_fisher_information}
  I\bigl( \mathcal{G}(a,p) \bigr) = \begin{cases}
                                      \frac{p^2}{a^2} \frac{\Gamma(2-1/p)}{\Gamma(1/p)}, & p > \frac{1}{2} \\
                                      \infty, & p \leq \frac{1}{2}
                                     \end{cases}
\end{equation}
denotes the Fisher Information of the zero-mean generalized normal distribution \cite[Chapter 3.2.1]{Kassam_2012}. 

Using the result in Section~\ref{sec:result}, an alternative lower bound can be obtained via a detour over the KL divergence between a Gaussian and a generalized Gaussian. The latter is given by \cite[eq.~(17)]{Do_2002}
\begin{align*}
  \label{eq:gg_kl}
  \DKL{\mathcal{G}(a,p)}{\mathcal{G}(a_0,2)} 
  &= \log\frac{p a_0}{2 a} \frac{\Gamma(1/2)}{\Gamma(1/p)}  + \frac{a^2}{a_0^2}\frac{\Gamma(3/p)}{\Gamma(1/p)} - \frac{1}{p}.
\end{align*}
In order to find the best Gaussian reference distribution, this distance needs to be minimized w.r.t.~$a_0$. It is not hard to show that this is accomplished by choosing  $a_0= \sqrt{\frac{2\Gamma(3/p)}{\Gamma(1/p)}} a$ (i.e., the KL divergence is minimized if the second moments agree), so that the following divergence can be defined:
\begin{align*}
  \DKL{\mathcal{G}(a,p)}{\mathcal{N}} &\coloneqq \min_{a_0 > 0} \; \DKL{\mathcal{G}(a,p)}{\mathcal{G}(a_0,2)} \\
  &= \log \frac{p}{\sqrt{2}} \sqrt{ \frac{\Gamma(3/p)}{\Gamma(1/p)}}\frac{\Gamma(1/2)}{\Gamma(1/p)}  + \frac{1}{2} - \frac{1}{p} \\
  &\eqqcolon d_{\text{KL}}(p),
\end{align*}
where $d_{\text{KL}} \colon \Reals_{> 0} \to \Reals_{> 0}$ is defined implicitly and only depends on $p$. Using this relation, bounds on the MMSE for any $X \sim \mathcal{G}(a,p)$ can be obtained by evaluating \eqref{eq:sigma_min} and \eqref{eq:sigma_max} at $\varepsilon = d_{\text{KL}}(p)$. 

\begin{figure}[t]
  \centering
  \includegraphics{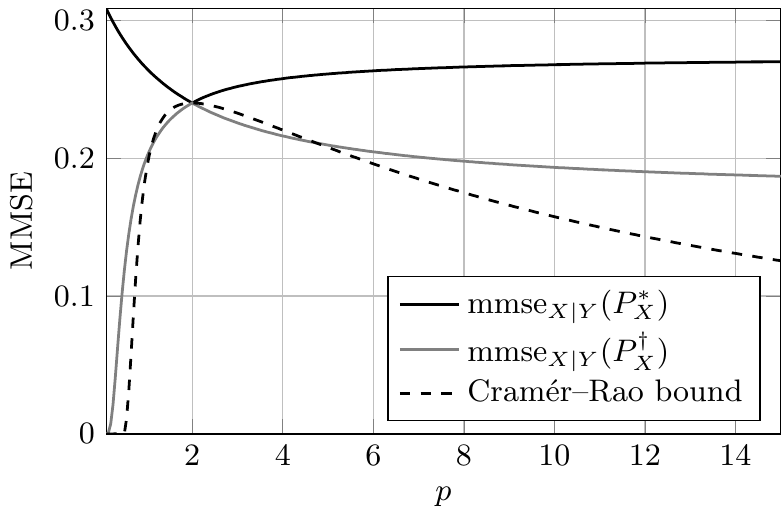}
  \caption{MMSE bounds vs.~the shape parameter $p$ of the generalized Gaussian distribution at SNR $=$ 5\,dB.}
  \label{fig:bounds_vs_p}
\end{figure}

An example of bounds on $\mmse_{X|Y}(\mathcal{G}(a,p))$ that were obtained via the KL divergence approach is shown in Fig.~\ref{fig:bounds_vs_p}. The SNR was set to 5\,dB and the Cram\'{e}r--Rao bound is plotted for comparison. Interestingly, neither of the bounds is uniformly tighter than the other. While the Cram\'{e}r--Rao bound is more accurate on the interval $\approx [1.1, 4.7]$, the bound in Section~\ref{sec:result} is tighter for all $p$ outside this interval. In particular for very large and very small values of $p$, the proposed bound is a significant improvement. Moreover, it can also be calculated for values $p \leq 0.5$, for which the Fisher Information in \eqref{eq:gg_fisher_information} is infinite so that the Cram\'{e}r--Rao bound becomes trivial. 


\bibliographystyle{IEEEbib}
\bibliography{refs}

\end{document}